# Partial Signatures and the Yoshida-Nicolaescu Theorem

J. C. C. Eidam, P. Piccione*

February 22, 2005

**Abstract**

In this article, we give a proof of the Yoshida-Nicolaescu Theorem by using the theory of partial signatures as in [8]. We do not impose the condition of non-degeneracy at the endpoints and use a natural definition of the Maslov index in this context. The proof here is more simple and direct and works in a more general context.

## 1 Introduction

Consider $X$ a closed oriented riemannian manifold *partitioned* by a hypersurface $Y$, that is, $X = X_+ \cup X_-$ where $X_\pm$ are manifolds with common boundary $Y$. If $P$ is a Dirac operator defined on a vector bundle $E$ over $X$, then the symbol of $P$ in the transversal variable at $Y$ gives a complex structure in the space $L^2(E|_Y)$. This complex structure induces a sympletic structure in $L^2(E|_Y)$ and the Cauchy data spaces $H_\pm(P)$ corresponding to the operator $P$ and the decomposition $X = X_+ \cup X_-$ form a Fredholm pair of lagrangian subspaces. These subspaces are, roughly speaking, formed by *restrictions* $u|_Y$ of solutions of the equation $Pu = 0$ in $X_\pm$.

In [15], T. Yoshida showed that if $\{P(s)\}_{s\in[0,1]}$ is a curve of Dirac operators in a three dimensional manifold such that $P(0)$ and $P(1)$ are invertible, then the spectral flow of the family $\{P(s)\}$ coincides with the Maslov index of the pair of lagrangian curves $(H_+(P(s)), H_-(P(s)))$. L. Nicolaescu extended this result to higher dimensions in [11], under the same hypotesis. In [5], M. Daniel proved this same formula in the case of degeneracy at the endpoints. To do this, he uses the Yoshida-Nicolaescu theorem in the non-degenerate case and a convention to compute the spectral flow and the Maslov index in the degenerate case. In this article, we give a simple and direct proof of the Yoshida-Nicolaescu theorem which works both in the degenerate and the non-degenerate case and avoids the use of conventions for calculation of the spectral flow and the Maslov index. We use the technique of partial signatures of [8], which has its roots in the work of M. Farber and J. Levine [7]. Since Dirac operators are unbounded, we use the definition of spectral flow for paths of such operators given by [2]; this important point was untouched in the other presentations of the problem. The sympletic structure related to the problem is strongly used, as in [6].

*Research supported by Fapesp, São Paulo, Brasil, grants 01/00046-3 and 02/02528-8. The second author is partially sponsored by Cnpq, Brasil.

## 2 Partial Signatures and Spectral Flow

We begin by describing the notion of spectral flow of a curve of (closed) Fredholm self-adjoint (un)bounded operators. Let $\{T(s)\}_{s\in[a,b]}$ be a curve of such operators, that is $T : [a, b] \to \mathscr{CF}^{\text{sa}}(\mathsf{H})$ is continuous, where $\mathscr{CF}^{\text{sa}}(\mathsf{H})$ denotes the space of (un)bounded Fredholm self-adjoint operators on the Hilbert space $\mathsf{H}$, endowed with the gap metric, as defined in [9], Chapter IV. Since a Fredholm operator cannot have $\lambda = 0$ in his essential spectrum and the discrete spectrum is continuous in the gap topology, there exists a partition $a = s_0 < \ldots < s_n = b$ and $\delta_1, \ldots, \delta_n > 0$ such that $\Sigma(T(s)) \cap [-\delta_j, \delta_j]$ is a finite set of eigenvalues of finite multiplicity, if $s \in [s_{j-1}, s_j]$, for $j = 1, \ldots, n$. If $m(s, \delta_j)$ is the number of eigenvalues of $T(s)$ in $[0, \delta_j]$, then the spectral flow of $\{T(s)\}_{s\in[a,b]}$ is defined as

$$\text{sf}(\{T(s)\}) = \sum_{j=1}^{n} \left\{ m(s_j, \delta_j) - m(s_{j-1}, \delta_j) \right\}.$$

The spectral flow is well-defined, i.e., does not depend on the particular partition and on the choice of the $\delta_j$'s. It is invariant by a homotopy which fixes the endpoints and additive by concatenation. More details can be found in [2].

Now, we introduce the theory of partial signatures developed in [8] to compute the spectral flow of curves of Fredholm self-adjoint operators. We begin with the bounded case.

Let $\{T(s)\}$ a $C^\infty$ family of bounded self-adjoint operators in $\mathsf{H}$, having a isolated degeneracy at $s = s_0$. A $C^\infty$ function $u : [s_0 - \varepsilon, s_0 + \varepsilon] \to \mathsf{H}$ is a *root function* for $u_0$ if $u(s_0) = u_0$. The order of the root function $u$ is the order of the zero of the function

$$[s_0 - \varepsilon, s_0 + \varepsilon] \ni s \mapsto T(s)u(s) \in \mathsf{H}$$

at $s = s_0$ and is denoted by $\text{ord}(u)$. We define for $k \geq 1$,

$$W_k(T, s_0) \doteq \{u_0 \in \ker T(s_0) : \text{ there exists a root function } u(s) \text{ for } u_0 \text{ such that } \text{ord}(u) \geq k\}. \tag{2.1}$$

The bilinear form

$$B_k(T, s_0)(u_0, v_0) \doteq \frac{1}{k!} \left\langle \frac{d^k}{ds^k}\bigg|_{s=s_0} T(s)u(s), v_0 \right\rangle, \tag{2.2}$$

is well-defined and hermitian for $u_0, v_0 \in W_k(T, s_0)$ and $u(s)$ a root function for $u_0$. The partial signatures of $\{T(s)\}$ at $s = s_0$ are defined as

$$n_k^-(T, s_0) \doteq n^-(B_k) , \ n_k^+(T, s_0) \doteq n^+(B_k) , \ \sigma_k(T, s_0) \doteq \sigma(B_k), \tag{2.3}$$

where $n^-(B_k)$, $n^+(B_k)$ and $\sigma(B_k)$ denote the *index*, *coindex* and the *signature* of the form $B_k$.[1]

If $\{T(s)\}$ is analytic and $\varepsilon > 0$ is sufficiently small, we can obtain real analytic functions $\lambda_1, \ldots, \lambda_N$, $N = \dim \ker T(s_0)$, defined in $[s_0 - \varepsilon, s_0 + \varepsilon]$ which represent the eingenvalues of $T(s)$ obtained by perturbation of the eingenvalue $\lambda = 0$ for $T(s_0)$ (see [9], Chapter VII-3). Now, it is easy to see that $n_k^+(T, s_0)$ ($n_k^-(T, s_0)$) is the number of $j$'s such that $\lambda_j(s_0) = \ldots = \lambda_j^{(k-1)}(s_0) = 0$ and $\lambda_j^{(k)}(s_0) > 0$ ($\lambda_j^{(k)}(s_0) < 0$). These observations lead to the following theorem.

---

[1]$n^+(B_k)$ ($n^-(B_k)$) is the number of positive (negative) eigenvalues of $B_k$ and $\sigma(B_k)$ is the difference $n^+(B_k) - n^-(B_k)$.

THEOREM 2.1 *Under the preceding hypotesis, we have*

$$\text{sf}(\{T(s)\}_{s\in[s_0-\varepsilon,s_0]}) = \sum_{k\geq 1} \left\{ n_{2k}^-(T,s_0) + n_{2k-1}^+(T,s_0) \right\},$$

$$\text{sf}(\{T(s)\}_{s\in[s_0,s_0+\varepsilon]}) = -\sum_{k\geq 1} n_k^-(T,s_0)$$

$$\text{sf}(\{T(s)\}_{s\in[s_0-\varepsilon,s_0+\varepsilon]}) = \sum_{k\geq 1} \sigma_{2k-1}(T,s_0).$$

The result above is valid in the more general case of a piecewise real-analytic curve, that is, a curve $\{T(s)\}_{s\in[s_0-\varepsilon,s_0+\varepsilon]}$ for which there exists a partition $s_0 - \varepsilon = x_0 < \ldots < x_n = s_0 + \varepsilon$ such that $\{T(s)\}_{s\in(x_j,x_{j-1})}$ is real analytic, for $j = 1, \ldots, n$.

The following proposition is useful for defining the Maslov index.

PROPOSITION 2.2 *Let $\{T(s)\}_{s\in[s_0-\varepsilon,s_0+\varepsilon]}$ be a $C^\infty$ family of bounded self-adjoint operators and $U : [s_0 - \varepsilon, s_0 + \varepsilon] \to \text{GL}(\mathsf{H})$ be a smooth curve. Assume that $U^*(s)T(s) = T(s)U(s)$ and $U(s_0)|_{\ker T(s_0)} = I$. Then the curve $\hat{T}(s) = T(s)U(s)$ has the same partial signatures in $s = s_0$ as $\{T(s)\}$; in particular, they have the same spectral flow.*

This theory has a important extension for families of elliptic operators which we describe now. Before getting to this point, let us recall some facts about the theory of holomorphic families of unbounded operators. The absolute reference for this is [9].

A one-parameter family $\{T(s)\}_{s\in I}$, $I \subset \mathbb{R}$, of closed densely defined operators on $\mathsf{H}$ is a holomorphic family of type (A) if all operators $T(s)$ have the same fixed dense domain $\mathsf{D}$ and the function $I \ni s \mapsto T(s)u \in \mathsf{H}$ is holomorphic for all $u \in \mathsf{D}$. The following lemma shows that under a simple hypothesis, it is easy to check that a family of operators is holomorphic of type (A).

LEMMA 2.3 *Let $\{T(s)\}_{s\in I}$ a one-parameter family of closed densely defined operators with common domain $\mathsf{D}$ and assume that there exists a norm $\|\cdot\|_\mathsf{D}$ on $\mathsf{D}$ such that all operators $T(s)$ are $\|\cdot\|_\mathsf{D}$-continuous. Then, $\{T(s)\}_{s\in I}$ is holomorphic of type (A) if and only if the map*

$$I \ni s \mapsto T(s) \in \mathscr{L}(\mathsf{D}, \mathsf{H}) \tag{2.4}$$

*is holomorphic, where $\mathscr{L}(\mathsf{D},\mathsf{H})$ denotes the Banach space of bounded linear operators $\mathsf{D} \to \mathsf{H}$ and $\mathsf{D}$ is endowed with the norm $\|\cdot\|_\mathsf{D}$.*

PROOF. Obviously, if 2.4 is holomorphic then $\{T(s)\}$ is holomorphic of type (A). If $\{T(s)\}$ is holomorphic of type (A), then, by definition, the map $I \ni s \mapsto T(s)u \in \mathsf{H}$ is holomorphic for all $u \in \mathsf{D}$. By the principle of uniform boundedness, we conclude that 2.4 is holomorphic. ∎

Now we are ready to formulate a extension of the theory of partial signatures to a case of our interest. Let $E$ be a complex riemannian vector bundle over a manifold and $\{P(s)\}_{s\in[0,1]}$ a family of first-order elliptic self-adjoint operators in $E$ which is a holomorphic family of type (A) in $L^2(E)$. Since all $P(s)$ can be regarded as bounded operators $H^1(E) \to L^2(E)$, lemma 2.3 shows that our hypothesis is equivalent to assume that the map $[0,1] \ni s \mapsto P(s) \in \mathscr{L}(H^1(E), L^2(E))$ is holomorphic.[2]

---

[2]$H^1(E)$ is the usual Sobolev space of order 1.

If $s = s_0$ is a isolated degeneracy of $\{P(s)\}$, a function $u : [s_0 − \varepsilon, s_0 + \varepsilon] \to L^2(E)$ is a *generalized root function*, or simply a *root function*, for $\{P(s)\}$ if $P(s_0)u(s_0) = 0$. The *order* of a root function $u$ is the order of the zero of the function

$$[s_0 − \varepsilon, s_0 + \varepsilon] \ni s \mapsto P(s)u(s) \in \mathscr{D}'(E)$$

at $s = s_0$ and is denoted by $\mathrm{ord}(u)$. Here, $\mathscr{D}'(E)$ denotes the space of distribution sections of $E$. Note that we do not assume that a root function takes values in $H^1(E)$, but in $L^2(E)$.

We define the spaces $W_k(P, s_0)$, $k \geq 1$, as in formula 2.1. We see that $u_0 \in W_k(P, s_0)$ if and only there exist root functon $u : [s_0 − \varepsilon, s_0 + \varepsilon] \to H^1(E)$ such that $u(s_0) = u_0$ and $\mathrm{ord}(u) \geq k$. In fact, by using the ellipticity of $P(s_0)$ and induction on $k$, we see that $u_0 \in W_k(P, s_0)$ if and only if there exists $u_0 = u_0^{(1)}, \ldots, u_0^{(k-1)} \in H^1(E)$ such that

$$\sum_{j=0}^{r} \binom{r}{j} P^{(j)}(s_0) u_0^{(r-j)} = 0,$$

for $r = 0, \ldots, k − 1$, where $P^{(j)}(s)$ denotes $j$-th derivative of the map $I \ni s \mapsto P(s) \in \mathscr{L}(H^1(E), L^2(E))$. In particular, $u(s) = \sum_{j=0}^{r-1} s^j u_0^{(j)} \in H^1(E)$ is a root function for $u_0$ of order $\geq k$. Thus, the bilinear form

$$B_k(P, s_0)(u_0, v_0) \doteq \frac{1}{k!} \left\langle \frac{d^k}{ds^k}\bigg|_{s=s_0} P(s)u(s), v_0 \right\rangle, \tag{2.5}$$

defined for $u_0, v_0 \in W_k(P, s_0)$ is well-defined and hermitian. Also in this case the partial signatures $n_k^\pm(P, s_0)$, $\sigma_k(P, s_0)$ of $\{P(s)\}$ at $s = s_0$ are defined by formulas 2.3.

Because $\{P(s)\}$ is analytic, if $\varepsilon > 0$ is sufficiently small we can obtain real analytic functions $\lambda_1, \ldots, \lambda_N$, $N = \dim \ker P(s_0)$, defined in $[s_0 − \varepsilon, s_0 + \varepsilon]$ which represent the eingenvalues of $P(s)$ obtained by perturbation of the eingenvalue $\lambda = 0$ for $P(s_0)$ (see [9]). Now, it is easy to see that $n_k^+(P, s_0)$ ($n_k^-(P, s_0)$) is the number of $j$'s such that $\lambda_j(s_0) = \ldots = \lambda_j^{(k-1)}(s_0) = 0$ and $\lambda_j^{(k)}(s_0) > 0$ ($\lambda_j^{(k)}(s_0) < 0$), and this implies the following theorem, analogous to theorem 2.1.

THEOREM 2.4 *Let* $\{P(s)\}_{s \in [s_0−\varepsilon, s_0+\varepsilon]}$ *a holomorphic family of type (A) of first-order elliptic self-adjoint operators acting on sections of E having a isolated degeneracy at* $s = s_0$. *Then*

$$\mathrm{sf}(\{P(s)\}_{s \in [s_0−\varepsilon, s_0]}) = \sum_{k \geq 1} \left\{ n_{2k}^-(P, s_0) + n_{2k-1}^+(P, s_0) \right\},$$

$$\mathrm{sf}(\{P(s)\}_{s \in [s_0, s_0+\varepsilon]}) = -\sum_{k \geq 1} n_k^-(P, s_0)$$

$$\mathrm{sf}(\{P(s)\}_{s \in [s_0−\varepsilon, s_0+\varepsilon]}) = \sum_{k \geq 1} \sigma_{2k-1}(P, s_0).$$

# 3 The Maslov Index

In this section, we define the Maslov index for a single path and for a pair of lagrangian paths in a infinite dimensional Hilbert space $\mathsf{H}$ endowed with a complex structure $J$ and the corresponding sympletic structure given by $\omega(\cdot, \cdot) = \langle J \cdot, \cdot \rangle$. A invertible linear operator preserving $\omega$ is said to be a sympletic isomorphism and the group of such operators is denoted by $\text{Sp}(H)$. The subgroup of $\text{Sp}(H)$ formed by the operators which commute with $J$ is denoted by $U(H_\mathbb{R})$.

A pair $(L_0, L_1)$ of lagrangians of $H$ is said to be *complementary* if $L_0 \cap L_1 = 0$ and $L_0 + L_1 = H$. The set of all lagrangians in $H$ is denoted by $\Lambda$ and the set of lagrangians complementary to $L_1$ is denoted by $\Lambda_0(L_1)$. Given $L_0 \in \Lambda$, $\mathscr{F}_{L_0}(\Lambda)$ is the set of lagrangians $L$ which form a Fredholm pair with $L_0$. It is an open set of $\Lambda$. The Maslov cycle with vertex at $L_0 \in \Lambda$ is the set of all $L \in \Lambda$ that intersects not trivially $L_0$ and is denoted by $\Sigma_{L_0}$.

DEFINITION 3.1 Let $(L_0, L_1)$ be complementary lagrangians. We define

$$\varphi_{L_0,L_1} : \Lambda_0(L_1) \to \mathscr{L}^{sa}(L_0)$$
$$L \mapsto P_{L_0} J S,$$

where $S : L_0 \to L_1$ has graph $\text{Gr}(S) = L$ and $\mathscr{L}^{sa}(L_0)$ denotes the Banach space of bounded self-adjoint operators in $L_0$.

It is easy to see that the set $\mathscr{A} = \{\varphi_{L_0,L_1}\}$, where $(L_0, L_1)$ varies in the set of all complementary pairs of lagrangians, forms a real analytic atlas for the space $\Lambda$. We can consider in $\Lambda$ the *gap* topology, that is, the topology obtained by means of the metric

$$\delta(L, L') = \|P_{L'}|_L\| = \sup_{\substack{u \in L \\ \|u\|=1}} \text{dist}(u, L'),$$

as in [9], Chapter IV-2. We note that the gap topology in $\Lambda$ coincides with the topology obtained by means of the atlas $\mathscr{A}$. In particular, $\Lambda$ is a Hausdorff space, and we conclude the following proposition.

PROPOSITION 3.2 *Endowed with the atlas $\mathscr{A}$, the space $\Lambda$ is a real analytic Banach manifold modeled in the space of self-adjoint operators in Hilbert space.*

We recall the method of [8] for defining the Maslov index. The *fundamental grupoid* of a topological space $X$ is the set of homotopy classes of curves $[\gamma]$ with fixed endpoints and is denoted by $\pi(X)$.

THEOREM 3.3 *Let $X$ be a topological space, $\{U_\alpha\}_{\alpha \in A}$ a open cover of $X$ and $G$ a group. Given a family $\{\psi_\alpha\}_{\alpha \in A}$ of homomorphisms $\psi_\alpha : \pi(U_\alpha) \to G$, there exists a unique homomorphism $\psi : \pi(X) \to G$ such that $\psi|_{\pi(U_\alpha)} = \psi_\alpha$ for all $\alpha \in A$ if and only if $\psi_\alpha|_{\pi(U_\alpha \cap U_\beta)} = \psi_\beta|_{\pi(U_\alpha \cap U_\beta)}$ for all $\alpha, \beta \in A$ such that $U_\alpha \cap U_\beta \neq \emptyset$.*

Our aim is define the Maslov index for curves in $\mathscr{F}_{L_0}(\Lambda)$ using theorem 3.3. We begin by putting

$$\psi_{L_0,L_1} : \pi(\Lambda_0(L_1) \cap \mathscr{F}_{L_0}(\Lambda)) \to \mathbb{Z}$$
$$\gamma \mapsto \text{sf}(\varphi_{L_0,L_1} \circ \gamma),$$

where sf denotes the spectral flow. Now we show that $\psi_{L_0,L_1}(\gamma) = \psi_{L_0,L'_1}(\gamma)$ if $\gamma : [0,1] \to \Lambda_0(L_1) \cap \Lambda_0(L'_1) \cap \mathscr{F}_{L_0}(\Lambda)$ is continuous. By homotopy invariance of the spectral flow, we can assume that $\gamma$ is real analytic by parts. Since

$$\varphi_{L_0,L'_1} \circ \varphi_{L_0,L_1}^{-1}(T) = T(I + P_{L'_1,L_0}^{L_1}(P_{L_0}J|_{L_1})^{-1}T)^{-1},$$

for all $T \in \mathscr{L}^{sa}(\mathsf{L}_0)$, we have

$$\varphi_{\mathsf{L}_0,\mathsf{L}_1'}(\gamma(s)) = \varphi_{\mathsf{L}_0,\mathsf{L}_1}(\gamma(s))U(s),$$

where $U(s) = (I + P^{\mathsf{L}_1}_{\mathsf{L}_1',\mathsf{L}_0}S(s))^{-1}$, $S(s) : \mathsf{L}_0 \to \mathsf{L}_1$ depends analytically on $s$ and $\mathrm{Gr}(S(s)) = \gamma(s)$. Since $\ker \varphi_{\mathsf{L}_0,\mathsf{L}_1}(\gamma(s)) = \gamma(s) \cap \mathsf{L}_0 = \ker S(s)$ and $U(s)|_{\ker S(s)} = I$, proposition 2.2 and theorem 3.3 given the following.

THEOREM 3.4  *For all $\mathsf{L}_1 \in \Lambda_0(\mathsf{L}_0)$, there exists a unique grupoid homomorphism $\mu_{\mathsf{L}_0}$ defined for curves in $\mathscr{F}_{\mathsf{L}_0}(\Lambda) \cap \Lambda_0(\mathsf{L}_1)$ such that*

$$\mu_{\mathsf{L}_0}(\gamma) = \mathrm{sf}(\varphi_{\mathsf{L}_0,\mathsf{L}_1} \circ \gamma).$$

*This homomorphism is called Maslov index.*

There is another important invariance property of the Maslov index described in next lemma.

LEMMA 3.5  *Given $\gamma : [0,1] \to \mathscr{F}_{\mathsf{L}_0}(\Lambda)$ continuous and $U : [0,1] \to \mathrm{Sp}(\mathsf{H})$ such that $U(s)(\mathsf{L}_0) = \mathsf{L}_0$, then $\mu_{\mathsf{L}_0}(U \cdot \gamma) = \mu_{\mathsf{L}_0}(\gamma)$, where $(U \cdot \gamma)(s) = U(s)(\gamma(s))$.*

PROOF.  First, assume $U$ is constant. By our next proposition, the group $\mathrm{Sp}(\mathsf{H}, \mathsf{L}_0) = \{U \in \mathrm{Sp}(\mathsf{H}) : U(\mathsf{L}_0) = \mathsf{L}_0\}$ is contractible (in particular, path-connected). So, there is a continuous curve $V : [0,1] \to \mathrm{Sp}(\mathsf{H}, \mathsf{L}_0)$ such that $V(0) = I$ and $V(1) = U$. Putting $H(s,t) = V(s)(\gamma(t))$, we have that $\dim(H(s,j) \cap \mathsf{L}_0) = \dim(\gamma(j) \cap \mathsf{L}_0)$ for all $s \in [0,1]$ and $j = 0, 1$. The result follows from homotopy invariance of the Maslov index.

In the general case, defining $K(s,t) = U(st)(\gamma(s))$, we have $\dim(K(s,j) \cap \mathsf{L}_0) = \dim(\gamma(j) \cap \mathsf{L}_0)$ for all $s \in [0,1]$ and $j = 0, 1$, and the result follows by the previous paragraph. ∎

PROPOSITION 3.6  *The group $\mathrm{Sp}(\mathsf{H}, \mathsf{L}_0)$ is contractible.*

PROOF.  A operator $T \in \mathrm{Sp}(\mathsf{H}; \mathsf{L}_0)$ has the form

$$\begin{pmatrix} A & AJS \\ 0 & -J(A^{-1})^*J \end{pmatrix}.$$

with respect to the decomposition $\mathsf{H} = \mathsf{L}_0 \oplus \mathsf{L}_0^\perp$, where $A \in \mathrm{GL}(\mathsf{L}_0)$ and $S \in \mathscr{L}^{sa}(\mathsf{L}_0^\perp)$. This implies that the map

$$\begin{aligned} \Phi : \mathrm{Sp}(\mathsf{H}, \mathsf{L}_0) &\to \mathrm{GL}(\mathsf{L}_0) \times \mathscr{L}^{sa}(\mathsf{L}_0^\perp) \\ T &\mapsto (A, S) \end{aligned}$$

is a diffeomorphism, and the result follows from the Kuiper's theorem [10]. ∎

Now we define the Maslov index for pairs of curves in $\Lambda$. Given curves $\gamma_0, \gamma_1 : [a,b] \to \Lambda$, we say that $(\gamma_0, \gamma_1)$ is a Fredholm pair of curves if for each $s \in [a,b]$, $(\gamma_0(s), \gamma_1(s))$ is a Fredholm pair of lagrangians.

Under the complex structure $\hat{J} = (J, -J)$ and the induced sympletic form, the Hilbert space $\mathsf{H} \oplus \mathsf{H}$ is a symplectic space. If $\mathsf{L}_0, \mathsf{L}_1$ are lagrangians in $\mathsf{H}$, it is easy to see that $\mathsf{L}_0 \oplus \mathsf{L}_1 \subset \mathsf{H} \oplus \mathsf{H}$ is lagrangian. The diagonal $\Delta$ is lagrangian in $\mathsf{H} \oplus \mathsf{H}$.

DEFINITION 3.7 *Given a Fredholm pair $\gamma_0, \gamma_1 : [a,b] \to \Lambda$ of continuous curves, we define the Maslov index of the pair $(\gamma_0, \gamma_1)$ as*

$$\mu(\gamma_0, \gamma_1) = \mu_\Delta(\gamma_0 \oplus \gamma_1).$$

Clearly, the Maslov index for pairs is homotopy invariant and additive by concatenation. The following proposition gives the relation between the Maslov index for a single curve and for a pair of curves. For $\mathsf{L}_0 \in \Lambda$, we denote by $P_{\mathsf{L}_0}$ the orthogonal projection over $\mathsf{L}_0$.

PROPOSITION 3.8 *Given $\mathsf{L}_0 \in \Lambda$ and a continuous curve $\gamma : [a,b] \to \mathscr{F}_{\mathsf{L}_0}(\Lambda)$ we have*

$$\mu_{\mathsf{L}_0}(\gamma) = \mu(\gamma, \mathsf{L}_0). \tag{3.1}$$

PROOF. First, we assume that $\gamma(s)$ is complementary to $\mathsf{L}_0^\perp$, for all $s \in [a,b]$. Given $\mathsf{L}_1 \in \Lambda_0(\mathsf{L}_0) \cap \Lambda_0(\mathsf{L}_0^\perp)$, we have that $(\mathsf{L}_0^\perp \oplus \mathsf{L}_1, \Delta)$ and $(\mathsf{L}_0^\perp \oplus \mathsf{L}_1, \gamma(s) \oplus \mathsf{L}_0)$ are complementary pairs of lagrangians. Let $S(s) : \mathsf{L}_0 \to \mathsf{L}_0^\perp$ bounded operators such that $\mathrm{Gr}(S(s)) = \gamma(s)$; defining $\tilde{S}(s) : \mathsf{H} \to \mathsf{L}_0$ by $\tilde{S}(s) = S(s)P_{\mathsf{L}_0} - P_{\mathsf{L}_0^\perp}$, we have $\mathrm{Gr}(\tilde{S}(s)) = \mathrm{Gr}(S(s)) = \gamma(s)$. Defining $T(s) = (\tilde{S}(s), -P_{\mathsf{L}_0, \mathsf{L}_1}) : \Delta \to \mathsf{L}_0^\perp \oplus \mathsf{L}_1$, where $P_{\mathsf{L}_0, \mathsf{L}_1}$ is the projection over $\mathsf{L}_1$ parallel to $\mathsf{L}_0$, we obtain that $\mathrm{Gr}(T(s)) = \gamma(s) \oplus \mathsf{L}_0$. So,

$$\varphi_{\Delta, \mathsf{L}_0^\perp \oplus \mathsf{L}_1}(\gamma(s) \oplus \mathsf{L}_0) = \frac{1}{2}(W(s), W(s)),$$

where $W(s) = JS(s)P_{\mathsf{L}_0} - J(P_{\mathsf{L}_0^\perp} - P_{\mathsf{L}_0, \mathsf{L}_1})$. Since $W(s)|_{\mathsf{L}_0} = JS(s)$ and $W(s)|_{\mathsf{L}_0^\perp} = -JP_{\mathsf{L}_1, \mathsf{L}_0}|_{\mathsf{L}_0^\perp}$, with respect to the decomposition $\mathsf{H} = \mathsf{L}_0 + \mathsf{L}_0^\perp$, $W(s)$ has the form

$$W(s) = \begin{pmatrix} JS(s) & 0 \\ 0 & -JP_{\mathsf{L}_0, \mathsf{L}_1}|_{\mathsf{L}_0^\perp} \end{pmatrix}.$$

By definition of spectral flow,

$$\mu(\gamma, \mathsf{L}_0) = \mu_\Delta(\gamma(s) \oplus \mathsf{L}_0) = \mathrm{sf}(W(s)) = \mathrm{sf}(JS(s)) = \mu_{\mathsf{L}_0}(\gamma),$$

as desired.

Now we discuss the general case. Since the numbers involved in equation 3.1 are additive by concatenation, we can assume that exists $\mathsf{L}_2 \in \Lambda_0(\mathsf{L}_0)$ complementary to $\gamma(s)$ for all $s \in [a,b]$ (see [13] for the existence of complementary lagrangians). Taking $U \in \mathrm{Sp}(\mathsf{H})$ such that $U(\mathsf{L}_0) = \mathsf{L}_0$ and $U(\mathsf{L}_2) = \mathsf{L}_0^\perp$, we have that $U(\gamma(s))$ is complementary to $\mathsf{L}_0^\perp$, for all $s \in [a,b]$, so

$$\mu(\gamma, \mathsf{L}_0) = \mu(U \cdot \gamma, \mathsf{L}_0) = \mu_{\mathsf{L}_0}(U \cdot \gamma) = \mu_{\mathsf{L}_0}(\gamma),$$

by lemma 3.5. ∎

Now we describe the idea of the Maslov index for pairs given by Nicolaescu in [11] and compare it with ours. Given a Fredholm pair of curves $\gamma_0, \gamma_1 : [a,b] \to \Lambda$ and $\mathsf{L}_0 \in \Lambda$, since the map

$$\begin{aligned} \mathrm{U}(\mathsf{H}_\mathbb{R}) &\to \Lambda \\ U &\mapsto U(\mathsf{L}_0) \end{aligned}$$

is a fibration, there exists a lifting for the curve $\gamma_1$, that is, there exists a curve $\eta :$ $[a,b] \to U(H_\mathbb{R})$ such that $\eta(s)(L_0) = \gamma_1(s)$ for all $s \in [a,b]$. The following proposition shows that the curve $(\eta^{-1} \cdot \gamma_0)(s) = \eta(s)^{-1}(\gamma_0(s))$ can be used to compute the Maslov index.

PROPOSITION 3.9 *We have $\mu(\gamma_0, \gamma_1) = \mu_{L_0}(\eta^{-1} \cdot \gamma_0)$. In particular, the integer $\mu_{L_0}(\eta^{-1} \cdot \gamma_0)$ does not depend on the lifting $\eta$.*

PROOF. Just use lemma 3.5 and proposition 3.8. ∎

Let us explain how to use partial signatures to calcule the Maslov index in the real analytic case. We begin by considering just one curve $\gamma : [s_0 - \varepsilon, s_0 + \varepsilon] \to \mathscr{F}_{L_0}(\Lambda)$ having a isolated degeneracy at $s = s_0$. If $\varepsilon > 0$ is sufficiently small, we can choose $L_1 \in \Lambda$ complementary to $L_0$ and $\gamma(s)$, for all $s \in [s_0 - \varepsilon, s_0 + \varepsilon]$. So the curve of bounded self-adjoint operators $\varphi_{L_0, L_1} \circ \gamma$ has a isolated degeneracy at $s = s_0$, and we can give the following definition.

DEFINITION 3.10 We define the $L_0$-partial signatures of $\gamma$ at $s = s_0$ as

$$n_k^+(\gamma, t_0; L_0) \doteq n_k^+(\varphi_{L_0,L_1} \circ \gamma, t_0) \quad , \quad n_k^-(\gamma, t_0; L_0) \doteq n_k^-(\varphi_{L_0,L_1} \circ \gamma, t_0)$$

$$\sigma_k(\gamma, t_0; L_0) \doteq n_k^+(\gamma, t_0; L_0) - n_k^-(\gamma, t_0; L_0).$$

A argument similar to that used to show the compatibility condition of the homomorphisms $\psi_{L_0, L_1}$ in theorem 3.4 shows that the $L_0$-partial signatures of $\gamma$ at $s = s_0$ do not depend of $L_1$, so they are well-defined. The following proposition follows from theorem 2.1.

THEOREM 3.11 *If $\gamma : [a,b] \to \mathscr{F}_{L_0}(\Lambda)$ is a real analytic curve not entirely cointained in the Maslov cycle $\Sigma_{L_0}$, then*

$$\mu_{L_0}(\gamma) = \sum_{\substack{\gamma(t) \in \Sigma_{L_0} \\ a < t < b}} \left[ \sum_{k \geq 1} \sigma_{2k-1}(\gamma, t_0; L_0) \right] +$$

$$+ \sum_{k \geq 1} \left[ n_{2k}^-(\gamma, a; L_0) + n_{2k-1}^+(\gamma, a; L_0) \right] - \sum_{k \geq 1} n_k^-(\gamma, b; L_0).$$

The same ideas can be applied to a pair $(\gamma_0, \gamma_1)$ of analytic curves in $\mathscr{F}_{L_0}(\Lambda)$ having a isolated intersection at $s = s_0$.

DEFINITION 3.12 The partial signatures of the pair $(\gamma_0, \gamma_1)$ at $s = s_0$ are defined as

$$n_k^+(\gamma_0, \gamma_1, t_0) \doteq n_k^+(\gamma_0 \oplus \gamma_1, t_0; \Delta) \quad , \quad n_k^-(\gamma_0, \gamma_1, t_0) \doteq n_k^-(\gamma_0 \oplus \gamma_1, t_0; \Delta)$$

$$\sigma_k(\gamma_0, \gamma_1, t_0) \doteq n_k^+(\gamma_0, \gamma_1, t_0) - n_k^-(\gamma_0, \gamma_1, t_0).$$

THEOREM 3.13 *If $\gamma_0, \gamma_1 : [a,b] \to \Lambda$ form a Fredholm pair of real analytic curves and $\gamma_0(s) \cap \gamma_1(s) = 0$ for some $s$, then*

$$\mu(\gamma_0, \gamma_1) = \sum_{\substack{a < t < b \\ \gamma_0(t) \cap \gamma_1(t) \neq 0}} \left[ \sum_{k \geq 1} \sigma_{2k-1}(\gamma_0, \gamma_1, t_0) \right] +$$

$$+ \sum_{k \geq 1} \left[ n_{2k}^-(\gamma_0, \gamma_1, a) + n_{2k-1}^+(\gamma_0, \gamma_1, a) \right] - \sum_{k \geq 1} n_k^-(\gamma_0, \gamma_1, b).$$

# 4 The Yoshida-Nicolaescu Theorem

Consider, as in the introduction, a closed oriented riemannian manifold $X = X_+ \cup X_-$ partitioned by a hypersurface $Y = X_+ \cap X_-$. Let $E$ be a riemannian vector bundle over $X$ and $\{P(s)\}_{s \in [0,1]}$ a family of first-order elliptic self-adjoint operators acting in sections of $E$. We make two hypotesis about the family $\{P(s)\}$:

(H1) The operators $P(s)$ satisfy the following *weak* unique continuation property: if $P(s)u = 0$ and $u = 0$ in a open set $V$, then $u = 0$ in every component intersecting $V$;

(H2) The operators $P(s)$ have *cilindrical form*, that is, in a bicollar of $Y$, they have the form
$$P(s) = G \cdot \left( \frac{\partial}{\partial t} + B(s) \right),$$
where $t$ is the transversal coordinate to $Y$, $G$ is a endomorphism of $E$ independent of $s, t$ such that $G^2 = -I$, $G^* = -G$ and $B(s)$ is a elliptic self-adjoint operator in $Y$ independing on $t$;

(H3) $\{P(s)\}_{s \in [0,1]}$ is a holomorphic family of type (A) (see lemma 2.3 and the comments below it).

An important property possessed by the operators $P(s)$ which is derived from the hypotesis (H1) and (H2) is that a solution of the equation $P(s)u = 0$ in $X_\pm$ is entirely determined by its trace over $Y$.

The Cauchy data spaces $H_\pm(s)$ are defined as the image of the Calderón projectors associated to the operators $P(s)$. They coincide with the spaces
$$\{u|_Y \,:\, P(s)u = 0 \text{ in } X_\pm \text{ and } u \in H^{1/2}(E|_{X_\pm})\},$$
where $u|_Y$ means the trace of $u$ over $Y$. For more details, see [1], [14].

The operator $G$ induces a complex structure in the space $L^2(E|_Y)$; this structure is very important for the study of boundary value problems for operators which satisfies the hypotesis (H1) and (H2). We reprove here a statement that appears already in [11], [3] with a proof not entirely satisfactory.

PROPOSITION 4.1 *If a first-order elliptic self-adjoint operator $P$ satisfies the hypotesis (H1) and (H2), then the Cauchy data spaces $H_+(P), H_-(P)$ form a Fredholm pair of lagrangians in $L^2(E|_Y)$.*

PROOF. We begin by showing that $H_+(P)$ is a lagrangian; the same proof works for $H_-(P)$. Given $g_0, g_1 \in H_+(P)$, there exists $u_0, u_1 \in H^{1/2}(E|_{X_+})$ such that $Pu_j = 0$ and $u_j|_Y = g_j$, for $j = 0, 1$. This implies, by the Green-Stokes formula (see [12], Cap.XVII), that
$$\int_Y Gg_0 \cdot g_1 \, dy = \int_{X_+} Pu_0 \cdot u_1 \, dx - \int_{X_+} u_0 \cdot Pu_1 = 0.$$

So, $G(H_+(P)) \subset H_+(P)^\perp$. To show equality, we observe that using the technique of invertible doubles (as in [4]), we obtain a compact manifold $\tilde{X}$ containing $X_+$ (and, obviously, partitioned by $Y$), a vector bundle $\tilde{E}$ over $\tilde{X}$ extending $E$ and a *invertible* first-order elliptic operator $\tilde{P}$ on it whose restriction to $X_+$ is $P$. The well-known theory of invertible elliptic operators over partitioned manifolds say us that $H_-(\tilde{P})\,+$ $H_+(\tilde{P}) = L^2(E)$ (see [4], [12]). But, by construction, we have $H_-(\tilde{P}) = G(H_+(P))$ and $H_+(\tilde{P}) = H_+(P)$; this implies that $H_+(P)$ is lagrangian.

Now, we show that $(H_+(P), H_-(P))$ is a Fredholm pair. It is well-known that the intersection $H_+(P) \cap H_-(P)$ is finite dimensional: in fact, if $g \in H_+(P) \cap H_-(P)$, then there exists $u_\pm \in H^{1/2}(E|_{X_\pm})$ such that $u_\pm|_Y = g$. By the Green-Stokes formula, we see that $u$ defined as $u_\pm$ on $X_\pm$ is a weak solution for $P$. By ellipticity, we conclude that $u \in C^\infty(E)$ and $Pu = 0$ in the classical sense. Again, by ellipticity, the space of such functions is finite dimensional. to conclude the proof, since the spaces $H_\pm(P)$ are lagrangian, it suffices prove that the sum $H_+(P) + H_-(P)$ is closed. Consider the Atiyah-Patodi-Singer projector $\pi_\geq(P_+)$ corresponding to the operator $P_+ = P|_{X_+}$, that is, $\pi_\geq(P_+)$ is the orthogonal projection over the non-negative subspace corresponding to the operator $B$ of the hypothesis (H2). The difference $\pi_\geq(P_+) - \Pi(P_+)$ is a compact operator, where $\Pi(P_+)$ is the Calderón projector corresponding to the operator $P_+$ (see [4]). Considering the operator $P_-$, which admits the form

$$P_- = (-G) \cdot \left( \frac{\partial}{\partial t'} - B(s) \right),$$

in a collar of $Y$ in $X_-$, where $t' = -t$, we have $\pi_\geq(P_-) = (I - \pi_\geq(P_+)) + P_{\ker B}$. This implies that $\Pi(P_+) + \Pi(P_-) = I + K$, with $K$ compact. Since the image of the Calderón projectors $\Pi(P_\pm)$ coincide with $H_\pm(P)$, the result follows from the following lemma. ∎

LEMMA 4.2 *If $P, Q$ are projections[3] in $\mathsf{H}$ such that $P + Q = I + K$, with $K$ compact, then $\operatorname{Im} P + \operatorname{Im} Q$ is closed and finite codimensional.*

PROOF.  Just note that $\operatorname{Im} P + \operatorname{Im} Q \supset \operatorname{Im}(P + Q) = \operatorname{Im}(I + K)$ and this last space is closed and finite codimensional, by the Fredholm alternative. ∎

By proposition 4.1, it makes sense talk about the Maslov index of the pair of curves $(H_+(s), H_-(s))$. This is a pair of *analytic* curves in the corresponding lagrangian-grassmannian. In fact, just observe that the corresponding Calderón projectors form a analytic curve of operators. Despite its non-orthogonality, they can be used for checking regularity of $(H_+(s), H_-(s))$ by the trick described in [4], Lemma 12.8.

The following theorem is a version of the Yoshida-Nicolaescu theorem without assuming neither non-degenaracy at endpoints nor that the operators involved are Dirac operators.

THEOREM 4.3 *Under the hypotesis (H1), (H2), (H3), we have*

$$\operatorname{sf}(\{P(s)\}) = \mu(\{H_+(s)\}, \{H_-(s)\}). \tag{4.1}$$

PROOF.  By analyticity, the eigenvalues of $\{P(s)\}$ are identically null or have isolated zeros. Since identically null eigenvalues makes no contribution for spectral flow/Maslov index, we can assume that there exists only one degeneracy $s = s_0$ of the curve $\{P(s)\}$. Let $\mathsf{L}$ be a lagrangian in $L^2(E|_Y) \oplus L^2(E|_Y)$ complementary to $\Delta$ and $H_+(s) \oplus H_-(s)$ and operators $(U_+(s), U_-(s)) : \Delta \to L$ such that $\text{Gr}((U_+(s), U_-(s))) = H_+(s) \oplus H_-(s)$, for $s \in [s_0 - \varepsilon, s_0 + \varepsilon]$ (see [13]). Since

$$\varphi_{\Delta,L}(H_+(s) \oplus H_-(s)) = P_\Delta \circ (G, -G) \circ (U_+(s), U_-(s))$$
$$= \left(\frac{1}{2}G(U_+(s) - U_-(s)), \frac{1}{2}G(U_+(s) - U_-(s))\right),$$

we have

$$\mu(\{H_+(s)\}, \{H_-(s)\}) = \text{sf}\left(\frac{1}{2}G(U_+(s) - U_-(s)), \frac{1}{2}G(U_+(s) - U_-(s))\right).$$

---

[3]By a *projection* we means just a bounded idempotent operator in $\mathsf{H}$.

So, we must prove that $\text{sf}(\{P(s)\}) = \text{sf}(\{G(U_+(s) - U_-(s))\})$. Consider the map

$$\Phi : \Delta \cap (H_+(s_0) \oplus H_-(s_0)) \to \ker P(s_0)$$
$$(g, g) \mapsto u,$$

where $u$ is the unique solution of $P(s_0)u = 0$ in $X$ such that $u|_Y = g$ and $W_k(P, s_0), B_k(P, s_0)$ and $W_k(H, s_0), B_k(H, s_0)$ the objects corresponding to the curves of operators $\{P(s)\}$ and $\{G(U_+(s) - U_-(s))\}$, respectively, defined in section 2. Since $\Phi$ is a isomorphism and $\ker G(U_+(s_0) - U_-(s_0)) = \Delta \cap (H_+(s_0) \oplus H_-(s_0))$, is suficient to show that $\Phi(W_k(H, s_0)) \subset W_k(P, s_0)$ and $\Phi^*(B_k(P, s_0)) = B_k(H, s_0)$.[4]

Given $(g_0, g_0), (h_0, h_0) \in W_k(H, s_0)$ and $(g(s), g(s))$ a root-function of order $\geq k$ for the curve of operators $G(U_+(s) - U_-(s))$ such that $g(s_0) = g_0$, by the very definition of $U_\pm(s)$, we have

$$g(s) + U_+(s)g(s) \doteq f_+(s) \in H_+(s)$$
$$g(s) + U_-(s)g(s) \doteq f_-(s) \in H_-(s).$$

So, there exists $F_\pm(s) \in C^\infty(E|_{X_\pm}) \cap H^{1/2}(E|_{X_\pm})$ such that

$$\begin{cases} P(s)F_\pm(s) = 0 \\ F_\pm(s)|_Y = f_\pm(s) \end{cases}.$$

Putting

$$F(s) = \begin{cases} F_+(s) & \text{in } X_+ \\ F_-(s) & \text{in } X_- \end{cases},$$

we have that $F(s) \in L^2(E)$ and $F(s_0) = u_0$, where $u_0 = \Phi(g_0, g_0)$. Obviously, $P(s)F(s)$ is a distribution supported in $Y$. Let us compute this distribution: given $w \in C^\infty(E|_Y)$ and $\tilde{w} \in C^\infty(E)$ such that $\tilde{w}|_Y = w$, we have

$$\langle P(s)F(s), w \rangle = \int_X F(s) \cdot P(s)\tilde{w}\, dy$$
$$= \int_{X_+} F_+(s) \cdot P(s)\tilde{w}\, dx + \int_{X_-} F_-(s) \cdot P(s)\tilde{w}\, dx$$
$$= \int_Y G(f_+(s) - f_-(s)) \cdot w\, dy,$$

by the Green-Stokes formula. This implies that

$$P(s)F(s) = \delta_Y \otimes G(f_+(s) - f_-(s)) = \delta_Y \otimes G(U_+(s) - U_-(s))g(s),$$

---

[4]We observe that, in general, $W_{k+1}(\cdot, s_0) = \ker B_k(\cdot, s_0)$, for all $k$, so, using induction, the equality $\Phi^*(B_k(P, s_0)) = B_k(H, s_0)$ implies that the inclusion $\Phi(W_k(H, s_0)) \subset W_k(P, s_0)$ is in fact an equality.

where $\delta_Y$ denotes the delta distribution supported in $Y$. This implies that $u_0 \in W_k(P, s_0)$, so $\Phi(W_k(H, s_0)) \subset W_k(P, s_0)$, for all $k \geq 1$. Defining $v_0 = \Phi(g_0, g_0)$ we have

$$
\begin{aligned}
B_k(H, s_0)((g_0, g_0), (h_0, h_0)) &= \frac{1}{k!} \int_Y \frac{d^k}{ds^k}\Big|_{s=s_0} G(U_+(s) - U_-(s))g(s) \cdot h_0 \, dy \\
&= \frac{1}{k!} \Big\langle \frac{d^k}{ds^k}\Big|_{s=s_0} P(s)F(s), h_0 \Big\rangle \\
&= B_k(P, s_0)(u_0, v_0) \\
&= \Phi^*(B_k(P, s_0))((g_0, g_0), (h_0, h_0)),
\end{aligned}
$$

as desired. ∎

Before extending the result of the previous theorem, we make a few comments about a simple approximation result.

Let $X$ a Banach space and $\varphi : \mathbb{R} \to X$ a (continuous) curve. Given $\chi : \mathbb{R} \to \mathbb{R}$ a $C^\infty$ function with compact support which is identically one in $[0, 1]$, define

$$\varphi_\alpha(s) = \sqrt{\frac{\alpha}{\pi}} \int_\mathbb{R} e^{-\alpha(s-t)^2} \chi(t) \varphi(t) dt,$$

for $\alpha > 0$. It is easy to see that $\varphi_\alpha$ is a real-analytic function and given $\varepsilon > 0$, $\|\varphi_\alpha(s) - \varphi(s)\| < \varepsilon$ for $s \in [0, 1]$ and sufficiently small $\alpha$. Defining

$$\psi_\alpha(s) = \varphi_\alpha(s) + (1 - s)(\varphi(0) - \varphi_\alpha(0)) + s(\varphi(1) - \varphi_\alpha(1)),$$

we obtain a real-analytic curve $\psi_\alpha$ such that $\psi_\alpha(0) = \varphi(0)$ and $\psi_\alpha(1) = \varphi(1)$. Taking $\alpha_0 > 0$ sufficiently small such that $\|\varphi_{\alpha_0}(s) - \varphi(s)\| < \varepsilon$ for all $s \in [0, 1]$, we conclude that $\|\psi_{\alpha_0}(s) - \varphi(s)\| < 2\varepsilon$ for all $s \in [0, 1]$. It follows that given any continuous curve in a Banach space, it is possible to obtain a real-analytic (uniform) approximation of it, with the same endpoints.

LEMMA 4.4 *Let $\varepsilon > 0$ and $\{\nabla^s\}_{s \in [0,1]}$ a continuous one-parameter family of connections in a riemannian vector bundle $E$, endowed with a Clifford product structure, over a compact manifold $X$ partitioned by a hypersurface $Y$. Then, there exists a analytic one-parameter family of connections $\{\tilde{\nabla}^s\}_{s \in [0,1]}$ which coincides with $\{\nabla^s\}_{s \in [0,1]}$ at the endpoints and such that $\|\sigma^s - \tilde{\sigma}^s\| < \varepsilon$, $s \in [0, 1]$, where $\sigma^s$ and $\tilde{\sigma}^s$ denote the symbols of the Dirac operators corresponding to $\nabla^s$ and $\tilde{\nabla}^s$, respectively.*

PROOF. The space of all connections in $E$ is a Banach space endowed with the norm

$$\|\nabla\| = \sum_{\alpha=1}^N \sum_{j,k=1}^n \sup_{x \in U_\alpha} |\nabla_{e_j^\alpha} e_k^\alpha(x)|,$$

where $\{e_j^\alpha\}_{j=1}^n$, $n = \dim E$, is a orthonormal frame defined over the relatively compact open set $U_\alpha$ and $\{U_\alpha\}_{\alpha=1}^N$ is an open cover of $X$. Since convergence in this norm implies convergence of the symbols of the correspondent Dirac operators, the result follows from our previous comments. ∎

Now, the classical version of the Yoshida-Nicolaescu theorem can be easily obtained.

THEOREM 4.5 *Let $\{P(s)\}_{s\in[0,1]}$ a continuous one-parameter family of Dirac operators over a closed riemannian manifold X partitioned by a hypersurface Y. If $H_{\pm}(s)$ is the corresponding one-parameter family of Cauchy data spaces, then formula 4.1 holds.*

PROOF. Using lemma 4.4, perturb the family $\{P(s)\}_{s\in[0,1]}$ to obtain a real analytic family $\{Q(s)\}_{s\in[0,1]}$ of Dirac operators which satisfy the hypotesis (H2). Since Dirac operators satisfy the hypotesis (H1) (see [1], [4]) and spectral flow/Maslov index are invariant under small perturbations, the result follows. ∎

Departamento de Matemática,
Instituto de Matemática e Estatística
Universidade de São Paulo, Brasil.
E-mail addresses: zieca@ime.usp.br, piccione@ime.usp.br